\documentclass[12pt]{extarticle}

\usepackage[margin=1.5in]{geometry}  
\usepackage{graphicx}              
\usepackage{amsmath}               
\usepackage{amsfonts}              
\usepackage{amsthm}                
\usepackage{framed}
 \usepackage{amssymb}
 \usepackage{siunitx}

\newtheorem{thm}{Theorem}

\begin{document}

\nocite{*}

\title{\bf On the Sum of a Prime and a Square-free Number}

\author{\textsc{Adrian W. Dudek} \\ 
Mathematical Sciences Institute \\
The Australian National University \\ 
\texttt{adrian.dudek@anu.edu.au}}
\date{}

\maketitle

\begin{abstract}
\noindent We prove that every integer greater than two may be written as the sum of a prime and a square-free number. 
\end{abstract}

\section{Introduction}

We say that a positive integer is square-free if it is not divisible by the square of any prime. It was first shown by Estermann \cite{estermann} in 1931 that every sufficiently large positive integer $n$ can be written as the sum of a prime and a square-free number. In particular, he proved that the number $T(n)$ of such representations satisfies the asymptotic formula
\begin{equation} \label{estermannformula}
T(n) \sim \frac{c n}{\log n} \prod_{ p | n} \bigg( 1+\frac{1}{p^2-p-1} \bigg),
\end{equation}
where 
\begin{equation} \label{artin}
c=\prod_p \bigg( 1 - \frac{1}{p(p-1)} \bigg) = 0.3739558\ldots
\end{equation}
is a product over all prime numbers which is known as Artin's constant (see Wrench's computation \cite{wrench} for more details). In 1935, Page \cite{page} improved Estermann's result by giving a bound for the order of the error term in (\ref{estermannformula}), using estimates for the error in the prime number theorem for arithmetic progressions. Mirsky \cite{mirsky} extended these results in 1949 to count representations of an integer as the sum of a prime and a $k$-free number, that is, a number which is not divisible by the $k$-th power of any prime. More recently, in 2005, Languasco \cite{languasco} treated the possibility of Siegel zeroes (see Davenport \cite{davenport} for some discussion on this) with more caution so as to provide better bounds on the error.

The objective of the present paper is to prove the following theorem which completes the result of Estermann.

\begin{thm} \label{main1}
Every integer greater than two is the sum of a prime and a square-free number. 
\end{thm}

We prove this theorem by working in the same manner as Estermann, though we employ explicit estimates on the error term for the prime number theorem in arithmetic progressions. Specifically, if we let
$$\theta(x;q,a) = \sum_{\substack{p \leq x \\ p \equiv a (\text{mod } q)}} \log p$$
where the sum is over primes $p$, we require estimates of the form
\begin{equation} \label{error}
\bigg| \theta(x;q,a) - \frac{x}{\varphi(q)} \bigg| < \frac{\epsilon x}{\varphi(q)}
\end{equation}
where $\epsilon>0$ is sufficiently small and $x$ and $q$ are suitably ranged. Good estimates of this type are available due to Ramar\'{e} and Rumely \cite{ramarerumely}, but are only provided for finitely many values of $q$. It turns out, however, that this is sufficient, as the Brun--Titchmarsh theorem is enough for the remaining cases.

It should be noted that the best known asymptotic result in this area is that of Chen \cite{chen}, who proved that every sufficiently large even number is the sum of a prime and another number which is the product of at most two primes. Chen's result is, at present, the closest one seems to be able to get towards a proof of the Goldbach conjecture, which famously asserts that every even integer greater than two can be written as the sum of two primes. It would be interesting to see if one could modify the explicit proof of this paper to consider, instead of square-free numbers, numbers which are products of at most $m$ primes (where $m$ is fixed). 

\section{The Proof}

\subsection{The Setup}

Let $n\geq10^{10}$ be a positive integer. We will prove Theorem 1 for this range before resorting to direct computation for the remaining cases. As convention would have it, we let $\mu: \mathbb{N} \rightarrow \{-1,0,1\}$ denote the M\"obius function, where $\mu(n)$ is zero if $n$ is not square-free; otherwise $\mu(n)=(-1)^{\omega(n)}$ where $\omega(n)$ denotes the number of distinct prime factors of $n$. As such, we have that $\mu(1)=1$. For a positive integer $n$, it can be shown that the sum
$$\mu_2 (n) = \sum_{a^2 | n} \mu(a)$$
is equal to $1$ if $n$ is square-free and zero otherwise. Thus, it follows that the expression
$$T(n):=\sum_{p \leq n} \mu_2(n-p)$$
counts the number of ways that $n$ may be expressed as the sum of a prime and a square-free number. We will employ logarithmic weights so as to use the known prime number estimates with more ease, and so we define
$$R(n) := \sum_{p \leq n} \mu_2(n-p) \log p.$$
We note that $n$ is the sum of a prime and a square-free number if and only if $R(n) >0$. As such, the majority of this paper is dedicated to finding a lower bound for $R(n)$. The expression for $R(n)$ can be rearranged so as to involve weighted sums over the prime numbers in arithmetic progressions:
\begin{eqnarray*}
R(n) & = & \sum_{p \leq n} \log p \sum_{a^2 | (n-p)} \mu(a) \\
& = & \sum_{a \leq n^{1/2}} \mu(a) \sum_{\substack{p \leq n \\ a^2 | (n-p)}} \log p \\
& = & \sum_{a \leq n^{1/2}} \mu(a) \theta(n; a^2, n).
\end{eqnarray*}
We will split the range of this sum into three parts, for we shall use a different technique to bound each of them. Note first that if $(a,n)>1$, then we have trivially that $\theta(n;a^2,n) \leq \log n$. Thus, we may write
\begin{equation} \label{bigdaddy}
R(n) > \Sigma_1 + \Sigma_2 + \Sigma_3 -n^{1/2} \log n
\end{equation}
where
\begin{eqnarray*}
\Sigma_1 & = & \sum_{\substack{a \leq 13 \\ (a,n)=1}} \mu(a) \theta(n; a^2, n), \\
\Sigma_2 & = & \sum_{\substack{13 < a \leq n^A \\ (a,n)=1}} \mu(a) \theta(n; a^2, n), \\
\Sigma_3 & = & \sum_{\substack{n^A < a \leq n^{1/2} \\ (a,n)=1}} \mu(a) \theta(n; a^2, n),
\end{eqnarray*}
and $A \in (0,1/2)$ is to be chosen later to optimise our result. We will use the estimates of Ramar\'{e} and Rumely \cite{ramarerumely} to bound $\Sigma_1$; this is the reason for the specific range of $a$ in this sum. We will then use the Brun--Titchmarsh theorem to bound $\Sigma_2$. Finally, $\Sigma_3$ will be bounded using trivial estimates. 

\subsection{Arithmetic estimates}

We first consider the sum
$$\Sigma_1 = \sum_{ \substack{a \leq 13 \\ (a,n)=1}} \mu(a) \theta(n; a^2, n).$$
Theorem 1 of Ramar\'{e} and Rumely \cite{ramarerumely} provides estimates of the form
\begin{equation}
\bigg| \theta(n;a^2,n) - \frac{n}{\varphi(a^2)} \bigg| < \frac{\epsilon_a n}{\varphi(a^2)}.
\end{equation}
In particular, by looking through the square moduli in Table 1 of their paper, we have values of $\epsilon_a$ for all $1 \leq a \leq 13$ which are valid for all $n\geq 10^{10}$. We therefore have trivially that
\begin{eqnarray} \label{big1}
\Sigma_1 & > & n \sum_{ \substack{a \leq 13 \\ (a,n)=1}} \bigg( \frac{\mu(a)}{\varphi(a^2)} - \epsilon_a \frac{\mu^2(a)}{\varphi(a^2)} \bigg) \nonumber \\
& > &  n \bigg( \sum_{(a,n)=1} \frac{\mu(a)}{\varphi(a^2)} - \sum_{\substack{a>13 \\ (a,n)=1}} \frac{\mu(a)}{\varphi(a^2)}- \sum_{a \leq 13} \frac{ \epsilon_a \mu^2(a)}{\varphi(a^2)} \bigg).
\end{eqnarray}
We wish to estimate the three sums in the above parentheses. We denote the leftmost sum by $S_{n}$, and note that we can bound it below by Artin's constant (\ref{artin}) \textit{viz.}
\begin{eqnarray*}
S_{n} := \sum_{(a,n)=1} \frac{\mu(a)}{\varphi(a^2)} & = & \prod_{p \nmid n} \bigg(1- \frac{1}{\varphi(p^2)} \bigg) \\
& = & \prod_{p \nmid n} \bigg(1- \frac{1}{p(p-1)} \bigg)\\
& \geq & \prod_{p} \bigg(1 - \frac{1}{p(p-1)}\bigg) = c.
\end{eqnarray*}
Wrench \cite{wrench} has computed this constant to high accuracy; it will suffice for the purpose of Theorem \ref{main1} to note that $S_n \geq c>0.373$. 


We will, for the moment, neglect the middle sum in (\ref{big1}), for it shall be considered jointly with a term in the estimation of $\Sigma_2$. Thus, in our estimation of $\Sigma_1$, it remains to manually compute the upper bound for the rightmost sum. This is a straightforward task which is done in reference to Table 1 of Ramare and Rumely's paper \cite{ramarerumely}. We get that
$$ \sum_{a \leq 13} \frac{\epsilon_a \mu^2(a)}{\varphi(a^2)} < 0.005.$$

We now bring $\Sigma_2$ into the fray; the explicit Brun--Titchmarsh theorem (see Montgomery and Vaughan \cite{montgomeryvaughan}) provides the bound
$$\theta(n; a^2,n) < 2 \bigg(\frac{\log n}{\log(n/a^2)} \bigg) \frac{n}{\varphi(a^2)}.$$
Clearly, in the range $13 < a \leq n^{A}$ we may bound 
$$\frac{\log n}{\log(n/a^2)} \leq \frac{1}{1-2A}$$
and so we have the estimate that
$$\theta(n; a^2,n) = \frac{n}{\varphi(a^2)} + \epsilon \bigg(\frac{1+2A}{1-2A} \bigg)\frac{n}{\varphi(a^2)}$$
where $|\epsilon| < 1$. We may then bound $\Sigma_2$ from below by
$$\Sigma_2 > n \sum_{\substack{13 < a \leq n^A \\ (a,n) = 1}} \frac{\mu(a)}{\varphi(a^2)} - n \bigg(\frac{1+2A}{1-2A}\bigg) \sum_{\substack{13 < a \leq n^A \\ (a,n)=1}} \frac{ \mu^2(a)}{\varphi(a^2)}.  $$
We can then add this to our estimate for $\Sigma_1$ to get
\begin{eqnarray} \label{sigma1and2}
\Sigma_1 + \Sigma_2 & > & n \bigg( S_{n} - 0.005 - \sum_{\substack{a> n^A \\ (a,n)=1}} \frac{\mu(a)}{\varphi(a^2)} - \bigg(\frac{1+2A}{1-2A}\bigg) \sum_{\substack{13 < a \leq n^A \\ (a,n)=1}} \frac{ \mu^2(a)}{\varphi(a^2)} \bigg) \nonumber \\
& > & n \bigg(S_n -0.005 - \bigg(\frac{1+2 A}{1- 2 A} \bigg) \sum_{a > 13} \frac{\mu^2 (a)}{\varphi(a^2)}\bigg).
\end{eqnarray}
We can estimate the sum in the above inequality by writing it as follows:
\begin{eqnarray*}
\sum_{a > 13} \frac{\mu^2 (a)}{\varphi(a^2)} =\sum_{a=1}^{\infty} \frac{\mu^2(a)}{\varphi(a^2)} - \sum_{a\leq 13} \frac{\mu^2(a)}{\varphi(a^2)}.
\end{eqnarray*}
The infinite sum is less than 1.95 (see Ramar\'{e} \cite{ramare} for example), and the finite sum can be computed by hand to see that the sum in (\ref{sigma1and2}) is bounded above by 0.086. Thus
\begin{equation} \label{sigma12}
\Sigma_1 + \Sigma_2 > n \bigg(S_n - 0.005 - \bigg(\frac{1+2A}{1-2A}\bigg)(0.086)\bigg).
\end{equation}

For $\Sigma_3$, we have trivially that
\begin{eqnarray} \label{sigma3bound}
|\Sigma_3| & < & \sum_{n^{A} < a \leq n^{1/2}} \theta(n;a^2,n) \nonumber  \\
& < & \sum_{n^{A} < a \leq n^{1/2}} \bigg(1+\frac{n}{a^2} \bigg) \log n \nonumber \\
& < & n^{1/2} \log n + n \log n \sum_{n^{A} < a \leq n^{1/2}} \frac{1}{a^2} \nonumber \\
& < & n^{1/2} \log n + n^{1-2A} \log n + n \log n \int_{n^{A}}^{n^{1/2}} t^{-2} dt \nonumber \\
& = & n^{1-2A} \log n + n^{1-A} \log n.
\end{eqnarray}

\subsection{A Lower Bound for $R(n)$}

We can now provide an explicit lower bound for
$$R(n) \geq \Sigma_1 + \Sigma_2 + \Sigma_3 -n^{1/2} \log n.$$
We combine our explicit estimates (\ref{sigma12}) and (\ref{sigma3bound}) and divide through by $n$ to get that
\begin{eqnarray*}
\frac{R(n)}{n} & > & S_{n}- 0.005-\bigg(\frac{1+2A}{1-2A}\bigg)(0.086) \\
& - &n^{-1/2} \log n - n^{-2A} \log n - n^{-A} \log n.
\end{eqnarray*}
For sufficiently small $A$ and large $n$, the right hand side will be positive. For any $n$, we have $S_n > 0.373$; it is a simple matter to choose $A=1/4$ and verify that the right hand side is positive for all $n \geq 10^{10}$. That is, Theorem \ref{main1} is true for all integers $n \geq 10^{10}$.

It so remains to prove this result for all integers in the range $3 \leq n <10^{10}$. If $n$ is even, we have the numerical verification by Oliviera e Silva, Herzog and Pardi \cite{silva} that all even integers up to $4 \cdot 10^{18}$ can be written as the sum of two primes. Thus, every even integer greater than two may be written as the sum of a prime and a square-free number. 

Thus, we need to check that every odd integer $3 \leq n \leq 10^{10}$ can be written as the sum of a prime and a square-free number. For each $n$, we subtract a prime number $p$ and check that the result is squarefree. This is a straightforward computation; we ran this on \textsc{Mathematica} and it took just under 3 days on a 2.6GHz laptop. The computation was eased somewhat by subtracting primes which were close in size to $n$.  

\subsection*{Acknowledgements}

The author wishes to thank Dr Timothy Trudgian for many helpful conversations regarding the present paper. 

\clearpage

\bibliographystyle{plain}

\bibliography{bibliov5}

\end{document}